\documentclass[12pt,reqno,a4paper,twoside]{article}

\usepackage{amsmath,amsthm,amstext,amscd,amssymb,euscript,mathrsfs}
\usepackage{epsfig}
\usepackage{color}
\usepackage{url}
\usepackage{comment}

\newcommand{\Z}{\mathbb Z}
\newcommand{\R}{\mathbb R}
\newcommand{\N}{\mathbb N}

\newcommand{\E}{\mathbb E}

\renewcommand{\phi}{\varphi}
\newcommand{\la}{\ensuremath{\Lambda}}
\newcommand{\si}{\ensuremath{\sigma}}

\newcommand{\bord}{\sf{b}}

\newcommand{\ising}{\mathrm{\scriptscriptstyle{Ising}}}
\newcommand{\mising}{\mu^{\scriptscriptstyle{M}}_\infty}

\def\1{\mathds{1}}

\newtheorem{theorem}{{\small T}{\scriptsize HEOREM}}[section]
\newtheorem{corollary}{{\bf{\small C}{\scriptsize OROLLARY}}}[section]
\newtheorem{proposition}{{\bf{\small P}{\scriptsize ROPOSITION}}}[section]
\newtheorem{lemma}{{\bf{\small L}{\scriptsize EMMA}}}[section]
\newtheorem{remark}{{\bf{\small R}{\scriptsize EMARK}}}[section]
\newtheorem{definition}{{\bf{\small D}{\scriptsize EFINITION}}}[section]

\renewenvironment{proof}[1]
{\noindent{{\bf{\small{ P}{\scriptsize ROOF}}}.}\hspace{0.1cm} #1} {$\;\qed$\newline}

\newcommand{\beq}{\begin{eqnarray}}
\newcommand{\eeq}{\end{eqnarray}}

\newcommand{\ba}{\begin{align*}}
\newcommand{\ea}{\end{align*}}

\newcommand{\be}{\begin{equation}}
\newcommand{\ee}{\end{equation}}

\newcommand{\bl}{\begin{lemma}}
\newcommand{\el}{\end{lemma}}

\newcommand{\br}{\begin{remark}}
\newcommand{\er}{\end{remark}}

\newcommand{\bt}{\begin{theorem}}
\newcommand{\et}{\end{theorem}}

\newcommand{\bd}{\begin{definition}}
\newcommand{\ed}{\end{definition}}

\newcommand{\bp}{\begin{proposition}}
\newcommand{\ep}{\end{proposition}}

\newcommand{\bc}{\begin{corollary}}
\newcommand{\ec}{\end{corollary}}

\newcommand{\bpr}{\begin{proof}}
\newcommand{\epr}{\end{proof}}

\newcommand{\bi}{\begin{itemize}}
\newcommand{\ei}{\end{itemize}}

\newcommand{\ben}{\begin{enumerate}}
\newcommand{\een}{\end{enumerate}}

\newcommand{\caA}{{\mathcal A}}

\newcommand{\ld}{{\ldots}}

\newcommand{\caI}{{\mathcal I}}

\newcommand{\caK}{{\mathcal K}}

\newcommand{\caN}{{\mathcal N}}

\newcommand{\caP}{{\mathcal P}}

\newcommand{\caR}{{\mathcal R}}

\begin{document}

\title{
Thermodynamic formalism and\\ large deviations for multiplication-invariant potentials
on lattice spin systems
}

\author{
Jean-Ren\'{e} Chazottes$^{\textup{{\tiny(a)}}}$ and Frank Redig$^{\textup{{\tiny(b)}}}$\\
{\small $^{\textup{(a)}}$
Centre de Physique Th\'eorique,}\\
{\small CNRS-\'Ecole polytechnique}\\
{\small 91128 Palaiseau Cedex, France}\\
{\small $^{\textup{(b)}}$
Delft Institute of Applied Mathematics,}\\
{\small Technische Universiteit Delft}\\
{\small Mekelweg 4, 2628 CD Delft, Nederland}\\
}

\maketitle

\begin{abstract}
We introduce the multiplicative Ising model and prove basic properties of its thermodynamic formalism
such as existence of pressure and entropies. We generalize to one-dimensional ``layer-unique'' Gibbs measures for which the same results can be obtained. For more general models associated to a $d$-dimensional multiplicative invariant potential, we prove a large deviation theorem in the uniqueness regime for averages of multiplicative shifts of general local functions. This thermodynamic formalism is motivated by the statistical properties of multiple ergodic averages.
\bigskip

\noindent

\end{abstract}




\section{Introduction}

In \cite{ccgr} we studied large deviations of multiple ergodic averages for Ising spins with a product distribution.
We also established a relation between the partition functions associated to multiple ergodic averages and partition functions of associated shift-invariant spin systems.
The dimension of the corresponding lattice spin system is related to the number of primes involved in the
multiple ergodic average.
E.g. $\sum_{i} \si_i\si_{3i}+\si_{i}\si_{9i}$ leads to a one-dimensional model
with interaction of range $2$, whereas $\sum_{i} \si_i\si_{3i}+\si_{i}\si_{2i}$
leads to a two-dimensional model with nearest neighbor interaction.
Just as in the standard Gibbs formalism, starting from large deviation properties of sums of shifts of a continuous function under a product measure, one is lead naturally (by Cram\'{e}r transformation, see e.g.\ \cite[chapter 2]{dz}) from product measures towards the set of Gibbs measures with shift-invariant interactions, and in that class one can prove again the large deviation principle for sums of shifts of a continuous function. It is therefore natural to extend the study
of large deviation properties for multiple ergodic averages under product measures to a class of measures which form the natural multiplication invariant analogue of shift-invariant Gibbs measures. \\
It is the aim of this paper to make some steps in that direction.
Just as for shift-invariant Gibbs measures in lattice spin systems, the one-dimensional context where there is uniqueness of Gibbs measures (with finite range or not too slowly decaying interaction)  we expect that the uniqueness transfers to the multiplication invariant context.
We start in this paper with the study of the multiplicative Ising model, which is the simplest case to start with after having dealt with product measures. We show that there is a unique Gibbs measure and study its thermodynamic formalism: entropy, pressure, and large deviation rate functions.
We show that under this measure, there is a large deviation principle for ergodic sums of so-called first-layer functions. Next, we generalize
this to the context of multiplication invariant potentials in dimension one, where the associated Gibbs measure is still unique and decomposes as a product on independent layers of Gibbs measures with a corresponding shift-invariant potential.
Finally we generalize to higher dimensional models such that on each layer we have uniqueness. This leads to a class
of so-called ``layer-unique'' Gibbs measures for which we have the multiplicative analogue of relative entropy density
and a corresponding large deviation principle for the multiplicative empirical measure.

\section{Some notations and definitions}

\subsection{Shift and multiplicative shift}

We consider lattice spin systems with Ising $\pm 1$ spins on the positive integers.
We denote by $\N$ the set of positive integers and let $\N_0=\N\cup \{0\}$.
We simply denote by $[M,N]$ the lattice intervals $\{M,M+1,\ldots, N\}$ for $M,N\in\N_0$ such that $M<N$.
Configurations which are elements of $\Omega=\{-1,1\}^{\N}$ will be denoted by $\si,\eta,\xi$. We also set
$\Omega_0=\{-1,1\}^{\N_0}$. We use the notation $\sigma_{[M,N]}$ for the restriction of $\sigma$ to the
lattice interval $[M,N]$ ($\sigma_{[M,N]}$ is thus an element of $\{-1,+1\}^{[M,N]}$).
The {\em shift} is defined, as usual, by
\[
(\theta_i (\si))_j=\si_{i+j}
\]
for $\si\in \Omega_0$ ($i,j\in\N_0$). This is the natural way the semigroup $(\N,+)$ acts on 
$\Omega_0$.\\
We introduce the {\em multiplicative shift}  by setting
\[
(T_i \si)_j = \si_{ij}
\]
for $\si\in \Omega$ ($i,j\in\N$). This is the action of the semigroup $(\N,\times)$ (which is generated by the prime numbers). Note that the shift and the multiplicative shift {\em do not} commute.

\subsection{Invariant measures}

It is not {\em a priori} clear that there exist probability measures which are invariant by the
multiplicative shift, apart from the trivial case of product measures. We shall see non trivial
examples in this paper.\\
Product measures are also invariant under the shift. A natural question is
whether they are the only ones. In the realm of probability measures with positive entropy with respect to the shift,
this is indeed the case \cite{mike}. \\
More generally, stochastic processes $(X_n)$ that are both stationary (in the sense that $(X_n)$ and $(X_{n+k})$
have the same marginals for all $k\in\N$) and such that $(X_n)$ and $(X_{rn})$ have the
same marginals for all $r\in\N$ are called ``strongly stationary''  and were introduced in the
context of ergodic Ramsey theory. Their structure is known and involves Bernoulli systems and rotations on nilmanifolds as building blocks \cite{frantz}.

\subsection{Standard Ising model}

The standard Ising model on the lattice interval $[0,N]$ with boundary conditions $\pm$ on the right and free on the left is the probability measure on $\{-1,1\}^{[0,N]}$ given by
\[
\mu^{\ising, \emptyset, \pm }_{N}(\si_{[0,N]})=
\frac{e^{-H^{\ising; \emptyset, \pm}_N (\si_{[0,N]})}}{Z^{\ising;\emptyset, \pm}_N}
\]
where $N\geq 1$ and where the Hamiltonian is given by
\[
H^{\ising;\emptyset, \pm}_N (\si_{[0,N]})=-\beta\left( \sum_{i=0}^{N-1} J\si_i\si_{i+1} + \sum_{i=0}^{N} h\si_i+\si_{N}(\pm 1)\right).
\]
The parameters of this Hamiltonian are $\beta$ (inverse temperature), $J$ (coupling strength) and $h$ (magnetic field).
Finally $Z^{\ising;\emptyset, \pm}_N$ is the partition function given by
\[
Z^{\ising;\emptyset, \pm}_N=
\sum_{\si_0,\ldots,\si_N=\pm 1} e^{-H^{\ising;\emptyset, \pm}_N (\si_{[0,N]})}.
\]
The measures $\mu^{\ising, \emptyset, \pm }_N$ have a unique (not depending on the right boundary condition) weak limit as $N\to\infty$ which we denote by $\mu^{\ising}_\infty$.
The standard Ising model corresponds to the potential (in the sense of \cite{geo})
\[
U(\{i,i+1\}, \si)=-J\beta \si_i\si_{i+1}, \ U(\{i\}, \si)=-\beta h \si_i
\]
which is {\em shift-invariant}, i.e.,
\[
U(A+i, \si)= U(A, \theta_i\si),\; \forall i\in\N_0.
\]
Notice however that because we consider the Ising model on $\Omega_0$ with free boundary condition on the left end, the corresponding $\mu^{\ising}_\infty$ need not be shift-invariant (this is the case only when $h=0$).

\section{The multiplicative Ising model}

We define what we call the ``multiplicative Ising model'' with parameters $\beta$ (inverse temperature), $J$ (coupling strength) and $h$ (magnetic field) as the lattice spin system on
$\Omega$ with {\em formal} Hamiltonian
\be\label{multis}
H(\si) = -\beta\left(\sum_{i\in \N} J \si_i\si_{2i} + h\sum_{i\in \N} \si_i\right).
\ee
This corresponds to the potential
\be\label{multpot}
U(\{i\}, \si) = -\beta h \si_i, \ U(\{i,2i\},\si) = -J\beta \si_i\si_{2i}
\ee
and $U(A,\si)=0$ elsewhere. This potential is {\em invariant by the multiplicative shift} in the sense that
\[
U(iA, \si)=U(A, T_i\si)
\]
for all $A\subset \N$, $\si\in \Omega$, $i\in\N $. We shall simply say that it is {\em multiplication invariant}.

The potential $U$ is of course non-shift invariant and long-range. The usual uniqueness criteria
for one-dimensional lattice spin systems do not apply, as well as the Dobrushin uniqueness criterion
(even for small $\beta$); see \cite{geo} for the statements of these criteria.
However we shall prove later on that uniqueness holds.

The Hamiltonian corresponding to \eqref{multis} in the lattice interval $[1,2N]$ with boundary condition $\eta$
is defined as
\[
H_{N}^{\eta}(\si_{[1,2N]})=-\beta\left( \sum_{i=1}^N J \si_i\si_{2i} +
\sum_{i=1}^{2N} h\si_i \pm
\sum_{i=N+1}^{2N}(\si_i \eta_{2i})\right). 
\]
Further $H_{N}^{\pm}(\si_{[1,2N]})$ stands for the Hamiltonian with plus or minus boundary conditions.
Finally let
\[
H_{N}^{\emptyset} (\si_{[1,2N]})= \sum_{i=1}^N (-J\beta)\si_i\si_{2i} +
\sum_{i=1}^{2N} (-h\beta)\si_i
\]
be the Hamiltonian with free boundary conditions. Finally we introduce the corresponding finite-volume
probability measures $\mu^\emptyset_N$ and $\mu^{\eta}_N$:
\be\label{mumu}
\mu^\emptyset_N(\si_{[1,2N]})=\frac{e^{-H_{N}^{\emptyset}(\si_{[1,2N]})}}{Z^\emptyset_N}
\quad,\quad
\mu^\eta_N(\si_{[1,2N]})=\frac{e^{-H_{N}^{\eta}(\si_{[1,2N]})}}{Z^\eta_N}
\ee
where
\be\label{zz}
Z^\emptyset_N=Z^\emptyset_N(\beta,h)=\sum_{\si_1,\ldots,\si_{2N}=\pm 1} e^{-H_{N}^{\emptyset}(\si_{[1,2N]})}
\ee
\be\label{zzeta}
Z^\eta_N=Z^\eta_N(\beta,h)=\sum_{\si_1,\ldots,\si_{2N}=\pm 1} e^{-H_{N}^{\eta}(\si_{[1,2N]})}.
\ee

\subsection{Layer spins}
Let us put $h=0$ from now on. As we sill see later, the case $h\not=0$ can be taken into account by a simple change of the {\em a priori} measure.

In \cite{ccgr} we introduced a natural and useful relabeling of $\si$
spins. More precisely, to a configuration $\si\in\Omega$ we associate a sequence $(\tau^r)$ of configurations
in $\Omega_0$, indexed by odd numbers $r$ defined by
\be\label{koko}
\tau^{r}_i = \si_{r2^i},\; i\in\N_0.
\ee
We call $r\in 2\N_0+1$ the {\em layer index} and $i$ the one-dimensional coordinate in the layer. We thus have
the following picture for this layer representation:
\begin{center}
\begin{tabular}{llllllll}
\;$\vdots$ & & & & & & & \;\;\;\;\;\;\;$\vdots$\\ \\
$\tau^7$ =    & $\si_7$ & $\si_{14}$ & $\si_{28}$ & $\si_{76}$ & $\si_{152}$ & $\ldots$ & (layer index 7)\\ \\
$\tau^5$ =  & $\si_5$ & $\si_{10}$ & $\si_{20}$ & $\si_{40}$ & $\si_{80}$ & $\ldots$ & (layer index 5) \\ \\
$\tau^3$ =    & $\si_3$ & $\si_6$ & $\si_{12}$ & $\si_{24}$ & $\si_{48}$ & $\ldots$ & (layer index 3) \\ \\
$\tau^1$ =   & $\si_1$ & $\si_2$ & $\si_4$ & $\si_8$ & $\si_{16}$ & $\ldots$ & (layer index 1)
\end{tabular}
\end{center}

\medskip

Then we can write
\be\label{relab}
\sum_{i=1}^N \si_i\si_{2i}= \sum_{\substack{1\leq k\leq N\\ k\; \text{odd}}}\ 
\sum_{i : k2^i\leq N}\tau^{k}_i\tau^k_{i+1}.
\ee
As a consequence, for $r\in 2\N_0+1$ given, under the free-boundary condition measure $\mu^\emptyset_N$, we have
\[
\psi_2(r/N):= \lfloor\log_2(N/r)\rfloor
\]
spins in layer $r$ which together form a {\em standard one-dimensional Ising model on the lattice interval
$[0, \psi_2 (r/N)]$}, with free boundary conditions at $0$ and at the right end.
Different layers are independent.

Adding plus or minus boundary conditions in \eqref{relab} yields
\beq\label{hamdec}
H^\pm_N (\si_{[1,2N]})&=& -\beta\left(\sum_{i=1}^N \si_i\si_{2i}\pm \sum_{i=N+1}^{2N} \si_i\right)
\nonumber\\
&=& 
-\beta \sum_{\substack{1\leq r\leq N\\  r\; \text{odd}}}\left(\left(\sum_{i=1}^{\psi_2 (r/N)}  \tau^r_i\tau^r_{i+1} \right) \pm \tau^r_{\psi_2(r/N)+1}\right)
\eeq
i.e., each term in the sum over $r$ gets exactly one extra term $\tau^r_{\psi_2(r/N)+1}$.
Consequently, for $k=2r-1$ given, we have once more $\psi_2(k/N)$ spins in layer $k$ which together
are a {\em standard one-dimensional Ising model on the lattice interval $[0, \psi_2 (k/N)]$}, with free boundary condition at $0$ but now with $\pm$ boundary condition at the right end.
As before, different layers are independent.

\subsection{Layer stationarity and multiplication invariance}

In this subsection we prove a general relation between layer-stationarity and multiplication invariance.
By theorem 1.1. in \cite{jen}, a non i.i.d. process which is multiplication-invariant
such a process cannot be stationary and ergodic under shifts. In our context, an example of such a multiplication invariant dependent measure is given by the multiplicative Ising model
with $h=0$, which indeed is not stationary under shifts, as we shall see below.

\bt\label{layerstatimpliesmultinv}
Let the relation between $\si$ and $\tau$ spins be as in \eqref{koko}. Suppose that
the $\{\tau^r, r\in 2\N_0+1\}$ form an i.i.d. sequence of stationary processes, i.e.,
for every $r\in 2\N_0+1$, $\tau^r=\{\tau^r_i:i\in \N_0\}$ is a stationary process, and for different
$r$'s, $\tau^r$ are independent. Then the distribution of the corresponding $\si$ is multiplication
invariant.
\et
\bpr
We have to show that under the conditions of the theorem, for every finite collection of numbers $p_1,\ldots, p_k\in \N$, and $m\in \N$ the joint distribution of
\[
(\si_{mp_1}, \ld, \si_{mp_k})
\]
coincides with that of
\[
(\si_{p_1}, \ld ,\si_{p_k}).
\]
Write $p_i= r_i 2^{v_i}$ with $r_i\in 2\N_0+1$, $v_i\in\N_0$, and $m=s2^u$. Then, using \eqref{koko} we have to prove that the joint
distribution of
\[
\tau^{sr_i}_{v_i+u},\  i\in \{1,\ld, k\}
\]
coincides with that of
\[
\tau^{r_i}_{v_i},\ i\in \{1,\ld, k\}.
\]
Denote $r_1=r_{n_1}< r_{n_2}<\cdots< r_{n_\ell}$ such that $\{r_1,\ld r_k\}=\{r_{n_1}, \ld, r_{n_\ell}\}$.
Further denote
\[
X^{w} =(\tau^{sr_{i}}_{v_i+u}:1\leq i\leq k, r_i=r_{n_w})
\]
and
\[
Y^{w} =(\tau^{r_{i}}_{v_i}:1\leq i\leq k, r_i=r_{n_w}).
\]
Then, by the independence of different layers, the joint distribution of
\[
\tau^{sr_i}_{v_i+u}, \ i\in \{1,\ld, k\}
\]
coincides with the joint distribution of
\[
\otimes_{w=1}^\ell X^w
\]
where $\otimes$ denotes independent joining.
Similarly, the joint distribution of
\[
\tau^{r_i}_{v_i},\  i\in \{1,\ld, k\}.
\]
coincides with that of
\[
\otimes_{w=1}^\ell Y^w.
\]
Therefore, it remains to show that for each $w$
the distributions of $X^w$ and $Y^w$ coincide, but this in turn follows from the assumptions of the theorem, which imply that the layers $sr_i$ and $r_i$ have the same stationary distribution.
\epr

\subsection{The infinite-volume limit $\mising$}

As a consequence of the correspondence between the $\si$ and the $\tau$ spins, and the existence of
the infinite-volume limit in each layer of $\tau$ spins we have the following.
\bt\label{layerthm}\leavevmode\\
\begin{enumerate}
\item Unique limit measure:
The measures $\mu_N^\eta$ have a unique ($\eta$-independent) weak limit (as $N\to\infty$)
denoted by $\mising$.  
This measure is called the multiplicative Ising measure on $\{-1,1\}^\N$. As a consequence, for the infinite-volume specification built to the potential \eqref{multpot} corresponds a unique infinite-volume consistent Gibbs measure
(in the sense of \cite{geo})
given by the same $\mising$.
\item Independent Ising layers:
Under $\mising$, the $\tau$-spins defined by
\[
\tau^{r}_i= \si_{r2^i}
\]
for $r$ odd, $i\in \N_0$, are independent and distributed according to the
standard Ising model measure $\mu^{\ising}_\infty$ with free boundary condition on the left.
\item Multiplication invariance:
The measure $\mising$ is multiplication invariant, i.e.,
for all $i\in \N$, $\si$ and $T_i\si$ have the same distribution.
\end{enumerate}
\et
\br
It is also easy to see that for $h=0$ the distribution of $\si_i$ does not depend on $i$ (single marginal stationarity) but e.g. the distribution of $\si_i\si_{i+1}$ does depend on $i$ (no full stationarity).
\er
The infinite-volume measure $\mu^{\ising}_\infty$ is a Markov measure, by the nearest neighbor character of the interaction. The transition matrix of the corresponding Markov chain is given by
\[
Q(a,b)=\frac{ e^{G(a,b)} \langle e_b, \tilde{e}\rangle}{\lambda \langle e_a,\tilde{e}\rangle}
\]
where 
\begin{itemize}
\item[-] $e_+$ is the unit vector $(1,0)$, $e_-$ the unit vector $(0,1)$;
\item[-] $\langle\cdot,\cdot\rangle$ denotes innerproduct;
\item[-] for $a, b\in \{-1, +1\}$
\[
G(a,b)= \beta(J  ab + h  b);
\]
\item[-] $\lambda>0$ denotes the maximal eigenvalue of the transfer matrix $K$ with elements given by
\[
K(a,b)= e^{G(a,b)}
\]
with corresponding eigenvector $\tilde{e}$.
\end{itemize}
The initial measure $\pi$ of this Markov chain is given by the distribution of $\si_1=\tau^1_0$, i.e., the first spin on
every layer:
\begin{eqnarray}
\nonumber
\pi(+1)&=&\mu^{\ising}_\infty (\si_0=+1)
\\
\nonumber
&=&
\lim_{N\to\infty}\frac{\sum_{\si_1, \si_N=\pm 1} e^{\beta h} e^{\beta J\si_1} e^{\beta h\si_1} K^{N-1}(\si_1,\si_N)}
{\sum_{\si_0, \si_N= \pm 1} e^{\beta h\si_0}  K^{N}(\si_0,\si_N)}
\\
&=&
\frac{e^{\beta h}\sum_{a,b=\pm 1}
e^{\beta J a} e^{\beta h a} \langle e_a, \tilde{e}\rangle\langle\tilde{e}, e_b\rangle}
{\sum_{a,b=\pm 1}\lambda  e^{\beta h a} \langle e_a, \tilde{e}\rangle\langle\tilde{e}, e_b\rangle}\cdot
\label{plouc}
\end{eqnarray}
For $h=0$ and using $\tilde{e}=\frac{1}{\sqrt{2}}(1,1)$, $\lambda= e^{-\beta J}+e^{\beta J}$, this gives $\pi(+1)=1/2$, which coincides with the stationary distribution of the Markov chain with transition matrix $Q$. In that case, the distribution on layers $\mu^{\ising}_\infty$ is stationary under the shift.
As a consequence, by theorem \ref{layerstatimpliesmultinv}, the measure $\mising$ is multiplication invariant for $h=0$.

This is no longer the case for $h\not= 0$. Notice that except for $J=0$, the measure $\mising$ is not stationary under the shift.
E.g. the joint distribution of $\si_1,\si_2$ and $\si_3,\si_4$ are not equal because $\si_1,\si_2$ are two neighboring spins on the same layer, whereas $\si_3,\si_4$ are on different layers and hence independent.

From the Markov property of $\mu^{\ising}_\infty$ we have the following formula for the cylinders of the layer Gibbs measures
\be\label{onedcil}
\log\mu^{\ising}_\infty (\eta_0, \ldots, \eta_k)= \log\mu(\eta_0) + \sum_{i=0}^{k-1} \log Q(\eta_i, \eta_{i+1})
\ee
with the convention that the sum is zero if empty.
This formula is useful, in e.g.\ the computation of the entropy of the multiplicative Ising model.

\subsection{Free energies}

Let $\bord\in \{\emptyset, \pm\}$. We are going to compute the free energies
\[
f^{\bord} = \lim_{N\to\infty} \frac1N \log Z^{\bord}_N 
\]
where the $Z^{\bord}_N$'s are defined in \eqref{zz} and \eqref{zzeta}.
Letting
\[
Z^{\ising;\emptyset, \bord}_N=
\sum_{\si_0,\ldots,\si_N=\pm 1} e^{-H^{\ising;\bord, \pm}_N (\si_{[0,N]})}
\]
we get, using \eqref{koko} and \eqref{hamdec},
\be\label{partfac}
Z_N^{\bord} = \prod_{\substack{k\leq N\\ k\ \text{odd}}} Z^{\ising; \emptyset,\bord}_{\psi_2 (k/N)}.
\ee

The following lemma will be useful now and at several places later.
\bl\label{KOROA}
Let $\phi:\N_0\to\R$ be a measurable function such that there exist $C>0$ and $q>0$ such that
$|\phi(n)|\leq C n^q$ for all $n\in\N_0$.
Then we have
\be\label{fformula}
\lim_{N\to\infty}\frac1N \sum_{\substack{1\leq i\leq N\\ i\;\textup{odd}}} \phi\left(\lfloor\log_2(N/i)\rfloor\right)
=\sum_{p=0}^\infty\frac{1}{2^{p+2}}\ \phi(p).
\ee
\el
\bpr
Since $|\phi(x)|<x^q$, it suffices to take the limit along the subsequence $N=2^K$. Then
\begin{align*}
& \lim_{N\to\infty}\frac1N \sum_{\substack{1\leq i\leq N\\ i\;\textup{odd}}} \phi\left(\lfloor\log_2(N/i)\rfloor\right)
= \lim_{K\to\infty}\frac1{2^K}\sum_{0\leq s \leq K-1} \sum_{\substack{r=2^s \\r\; \text{odd}}}^{2^{s+1}-1} \phi(K-s-1)\\
&=\lim_{K\to\infty}\frac1{2^K}\sum_{1\leq s \leq K-1}  \phi(K-s-1) 2^{s-1}
= \lim_{K\to\infty} \sum_{1\leq s\leq K-1}\frac{1}{2^{K-s+1}}\ \phi(K-s-1)\\
&=\lim_{K\to\infty}\sum_{p=0}^{K-2}\frac{1}{2^{p+2}}\  \phi(p)
= \sum_{p=0}^{\infty}\frac{1}{2^{p+2}} \ \phi(p).
\end{align*}
\epr

As an application of lemma \ref{KOROA} and \eqref{partfac}, we obtain for the free energies of the multiplicative Ising model
\be\label{freefree}
f^{\bord}= \sum_{p=0}^\infty\frac{1}{2^{p+2}} \log Z^{\ising;\emptyset,\bord}_p
\ee
with $\bord\in \{\emptyset, \pm\}$.
We could derive more explicit expressions for $f^{\bord}$ (see \cite{ccgr} for instance), but here we just notice that, contrary to what one is used to
in the shift-invariant context, here the free energy is {\em depending on the boundary conditions}.
This is due to the non-shift invariant and long-range character of the interaction.

\section{Large deviation properties of $\mising$}

Gibbs measures with shift-invariant potentials satisfy nice large deviation properties, where the large deviation rate function is given by the relative entropy density, and the corresponding logarithmic moment-generating function given by a difference of free energies, see e.g. \cite{geo}.
In the present context, the natural invariance is multiplicative rather than additive, and so other large deviation properties will appear, and the natural quantities that are satisfying large deviation properties will be finite-volume Hamiltonians associated to a multiplicatively invariant potential.

\subsection{Free boundary conditions}

As a warming-up example we consider the large deviations of the normalized sums
\[
\frac{S^{(2)}_N}{N}= \frac1N\sum_{i=1}^N \si_i\si_{2i}
\]
under the measures $\mu^\emptyset_N$ on $\{-1,1\}^{[1,2N]}$ defined in \eqref{mumu}, each of them corresponding to the multiplicative Ising model in the lattice interval $[1,2N]$ with free boundary conditions.
In the shift-invariant context the large deviation rate functions, as well as the entropy in the thermodynamic limit, do not
depend on the boundary conditions. Here this is not the case. The free-boundary case is the easiest.

The free energy partition function is related to the free energy partition function of the standard Ising model via the correspondence
\eqref{koko}.

We can use the G\"{a}rtner-Ellis theorem (see \cite{dz}) and first compute, using \eqref{freefree}
and \eqref{zz}
\beq\label{freeendif}
F(t)&: =&\lim_{N\to\infty}\frac1N \log\E_{\mu^\emptyset_N}\left(e^{tS^{(2)}_N}\right)
\nonumber\\
&=&
\lim_{N\to\infty}\frac1N \log\frac{Z^\emptyset_N(\beta+t)}{Z^\emptyset_N(\beta)}
=
f^\emptyset(\beta+t)-f^\emptyset(\beta)
\nonumber\\
&=&
\sum_{p=0}^\infty\frac{1}{2^{p+2}}
\log \frac{Z^{\ising;\emptyset,\emptyset}_p(\beta+t)}{Z^{\ising;\emptyset,\emptyset}_p(\beta)}
\eeq
where $Z^{\ising;\emptyset,\emptyset}_p(\beta)$ is the partition function of the standard Ising model
in $[0,p]$ with free boundary conditions:
\[
Z^{\ising,\emptyset,\emptyset}_p(\beta)= \sum_{\si_0,\ldots, \si_{p}=\pm 1} e^{\beta\sum_{i=0}^{p-1} \si_i \si_{i+1}}
\]
where $p\geq 1$.
From \eqref{freeendif} we obtain existence and differentiability of $F(t)$ and we thus can conclude that $S^{(2)}_N/N$ satisfies the large deviation principle under the measures $\mu^\emptyset_N$ with rate function $I$ given by the Legendre transform of $F$: $I(x)= \sup_{t\in\R} (tx-F(t))$.

Similarly, we can easily obtain a formula for the ``free boundary condition'' entropy in the thermodynamic limit, using
lemma \ref{KOROA}:
\begin{align*}
s^\emptyset(\beta) &: =
\lim_{N\to\infty}\frac1N \; s(\mu_N^\emptyset)
\\
&= -\lim_{N\to\infty}\frac1N\sum_{\si_1,\ldots, \si_{2N}=\pm 1} \mu^\emptyset_N(\si_{[1,2N]})\log \mu^\emptyset_N(\si_{[1,2N]})
\\
&=
-\lim_{N\to\infty}\frac1N\left( \beta\frac{d}{d\beta} \log Z^\emptyset_N -\log Z^\emptyset_N\right)\\
&=
\sum_{p=0}^\infty\frac{1}{2^{p+2}} \ s^{\ising;\emptyset}_{p+1}(\beta)
\end{align*}
where
\[
s^{\ising;\emptyset}_{p+1}(\beta):=\left( \beta\frac{d}{d\beta} \log Z^{\ising;\emptyset,\emptyset}_p(\beta) -
\log Z^{\ising;\emptyset,\emptyset}_p(\beta)\right)
\]
is the entropy of the standard Ising model with free boundary conditions in the lattice interval $[0,p]$.

\subsection{Kolmogorov-Sinai entropy}\label{KS}

In the context of shift-invariant Gibbs measures, the free energy does not depend on the boundary condition, and
therefore, neither do large properties of sums of shifts under a shift-invariant Gibbs measure. In the multiplication invariant context, the thermodynamic formalism is different,
as we have already witnessed in the computation of the free energies for $\pm$ boundary conditions, which depend on the boundary condition.

To start with the study of the thermodynamic properties of $\mising$,
let us first consider its Kolmogorov-Sinai (KS) entropy.
First notice that by the fact that $\mising$ factorizes over different layers of $\tau$ spins, we have
\[
\log \mising (\si_{[1,N]})=\sum_{\substack{1\leq r\leq N\\  r \ \text{odd}}}\log\mu^{\ising}_\infty
(\tau^r_1,\ldots, \tau^r_{\psi_2 (r/N)}).
\]
Denote by
\[
s^{\ising}_{k+1}=-\E_{\mu^{\ising}_\infty}\left(\log\mu^{\ising}_\infty(\tau_0,\ldots,\tau_k)\right)
\]
the entropy of cylinders of length $k+1$ under the measure $\mu^{\ising}_\infty$.
Then, using lemma \ref{KOROA}, we obtain the following explicit formula for the KS entropy of $\mising$:
\[
-\lim_{N\to\infty}\frac1N\; \E_{\mising}\log \mising (\si_{[1,N]})=
\sum_{k=0}^\infty \frac{1}{2^{k+2}}\ s^{\ising}_{k+1}.
\]
To obtain a more explicit formula in terms of the transfer matrix, we use the Markov structure of the layer Gibbs measure
$\mu^{\ising}_\infty$.

Now, using lemma \ref{KOROA} and  \eqref{onedcil} we obtain
\begin{align*}
& -s(\mising) =\lim_{N\to\infty}\frac1N\; \E_{\mising} (\log\mising(\si_{[1,N]}))\\
&=
\frac12 \; \E_{\mu^{\ising}_\infty}\left(\log \mu^{\ising}_\infty(\tau_0)\right)
+
\lim_{N\to\infty}\frac1N\sum_{\substack{1\leq r\leq N\\ r \; \text{odd}}}
\; \sum_{0\leq i\leq \psi_2(r/N)} \E_{\mu^{\ising}_\infty}(\log Q(\tau_i, \tau_{i+1}))
\\
&=
\frac12 \; \E_{\mu^{\ising}_\infty}\left(\log \mu^{\ising}_\infty(\tau_0)\right)
+ \sum_{k=1}^\infty \frac{1}{2^{k+2}}\ \sum_{i=0}^{k-1} \E_{\mu^{\ising}_\infty}(\log Q(\tau_i, \tau_{i+1})).
\end{align*}
Using the elementary formula
\[
\sum_{k=1}^\infty \frac{1}{2^{k+2}} \sum_{i=0}^{k-1} z_i= \sum_{i=0}^\infty \frac{z_i}{2^{i+2}}
\]
and $\pi(\tau_0)=\mu^{\ising}_\infty(\tau_0)$, where $\pi$ is defined in \eqref{plouc}, we obtain
the following explicit formula
\begin{align}
\nonumber
& s(\mising)  =  -\frac12 \big(\pi(+1)\log\pi(+1)+\pi(-1)\log\pi(-1)\big)\\
& \qquad \qquad - \frac12\sum_{a,b,c=\pm 1} \pi(a)  \caR(a,b) Q(b,c)\log(Q(b,c))
\label{ksentropyMISING}
\end{align}
where
\[
\caR(a,b):= \sum_{k=0}^\infty \frac1{2^{k+1}} \; Q(a,b)^k=\frac12 \left(1-\frac{Q(a,b)}{2}\right)^{-1}.
\]
\br
\ben
\item When $J=0, h=0$, the measure $\mising$ is nothing but the product measure
giving weight $1/2$ to $\pm 1$ whose entropy is $\log 2$. 
Formula \eqref{ksentropyMISING} indeed gives $\log 2$.
\item If $h=0$ then 
\[
s(\mising)= \frac12 + \frac12 H(\alpha)
\] 
where $H(\alpha):=-\alpha \log\alpha -(1-\alpha)\log (1-\alpha)$ and $\alpha:=(1+e^{-2\beta J})^{-1}$.
If we fix $\beta J$ such that $\E_{\mising}(\si_1\si_2)=\gamma$ then choosing
$\alpha=\frac{1-\gamma}{2}$ yields
\[
s(\mising)= \frac12 + \frac12 H\left(\frac{1-\gamma}{2}\right).
\]
This corresponds to the Hausdorff dimension of the level set (see \cite{fanfan})
\[
\left\{\si : \lim_{n\to\infty} \frac1n \sum_{i=1}^{n} \si_i\si_{2i}=\gamma\right\}
\]
and shows that $\mising$ is the natural measure concentrated on this level set. 
In general, however, the telescopic measures constructed in \cite{fanfan} tailored for the multifractal
analysis of multiple ergodic averages are different from the Gibbs measures we introduce in this paper
for the purpose of large deviations. 
\een
\er

\subsection{The pressure of first-layer functions}

In order to obtain large deviation results for ``ergodic averages'' of the form
$\frac1N\sum_{i=1}^N T_i f$ under the measure
$\mising$, i.e., the infinite-volume multiplicative Ising model, we will heavily rely on the
independent layer decomposition.
We will therefore have to restrict to functions
$f$ such that their ergodic sums $\sum_{i=1}^N T_i f$ are consistent with this independent layer structure.
\bd
A continuous function $f:\Omega\to\R$ of the $\si$'s is called a {\em first-layer function} if there exists a continuous $f^*:\Omega_0\to\R$ of the $\tau$'s such that $f(\si)= f^*(\tau^1)$, i.e., if $f$ depends only on the spins in the first layer.
\ed

For such a first-layer function we have
\beq\label{BARAK}
\sum_{i=1}^N T_i f(\si)&=& \sum_{\substack{1\leq k\leq N\\ k\; \text{odd}}}
\;\sum_{0\leq i\leq \psi_2(k/N)} T_{k2^i} f(\si)
\nonumber\\
&=&
\sum_{\substack{1\leq k\leq N\\ k\; \text{odd}}}\; \sum_{0\leq i\leq\psi_2(k/N)} f^* (\theta_i(\tau^k)).
\eeq

We then define the pressure of a first-layer function w.r.t. multiplicative Ising model as follows.
Let us first define
\[
P^k_{\mu^{\ising}_\infty} (f^*)= \log\E_{\mu^{\ising}_\infty}\left(e^{\sum_{i=0}^{k} \theta_i f^*}\right).
\]
Since ${\mu^{\ising}_\infty}$ is a one-dimensional Gibbs measure with nearest-neighbor interaction, we have the estimate
\be\label{pk}
P^k_{\mu^{\ising}_\infty} (f^*)\leq C k.
\ee
Next define
\[
\caP^{\scriptscriptstyle{M}} (f|\mising)=\lim_{N\to\infty}\frac1N\log\E_{\mising}\left( e^{\sum_{i=1}^N T_i f}\right).
\]
Then, using \eqref{BARAK}, \eqref{pk} and lemma \ref{KOROA} for $f$ a first-layer function we have,
using the independence of the $\tau$ spins for different layers, and the fact that they are distributed according to the one-dimensional Gibbs measure $\mu^{\ising}_\infty$
\[
\caP^{\scriptscriptstyle{M}} (f|\mising)= \sum_{k=0}^\infty \frac1{2^{k+2}}\; P^k_{\mu^{\ising}_\infty} (f^*).
\]
We have the following result.
\bt\label{ldthm1}
\leavevmode\\
\begin{itemize}
\item[\textup{(a)}]
Under the measure $\mising$, for every first-layer function $f$, the random variables
\[
X_N(f):=\frac1N\sum_{i=1}^N T_i f
\]
satisfy the large deviation principle
with rate function
\[
I_f(x)= \sup_{t\in\R}\left(tx-\caP^{\scriptscriptstyle{M}}(t f|\mising)\right).
\]
\item[\textup{(b)}]
Moreover, they satisfy the central limit theorem, i.e.,
\[
\frac{1}{\sqrt{N}}( X_N(f)-\E_\mu(X_N(f)))\to \caN (0,\si^2)
\]
where $\to$ means here convergence in distribution, and where
\[
\si^2=\frac{\textup{d}^2\caP^{\scriptscriptstyle{M}}(t f|\mu)}{\textup{d}t^2}\Big|_{t=0}=
\lim_{N\to\infty} \frac1N\sum_{i,j=1}^N\big(\E_{\mising} (T_i f T_j f)-(\E_{\mising} (f))^2\big).
\]
\end{itemize}
\et

\bpr
\begin{itemize}
\item[\textup{(a)}]
We have 
\[
\left| \frac{\textup{d}P^k_{\mu^{\ising}_\infty} (f^*)}{\textup{d}t}\right| \leq c\; k
\]
for some constant $c>0$ and for all $k$.
Hence the function $t\mapsto \caP^{\scriptscriptstyle{M}}(t f|\mising)$, $t\in\R$, is continuously differentiable. The result
follows from G\"artner-Ellis theorem \cite{dz}.
\item[\textup{(b)}]
By complete analyticity of one-dimensional lattice models with finite-range interaction \cite{dobshlos}, it follows
that there exists a neighborhood $\mathcal{V}\subset \mathbb{C}$ of the origin and $C>0$ such that for all $k\in\N$
\[
\sup_{z\in \mathcal{V}} \left| \frac1k P^k_{\mu^{\ising}_\infty} (zf^*)\right|\leq C.
\]
Therefore the map $z\mapsto \caP^{\scriptscriptstyle{M}} (zf|\mising)$ is well defined for $z\in\mathcal{V}$ and one can apply
Bryc's theorem \cite{bryc}.
\end{itemize}
\epr

We can push the large deviation result of theorem \ref{ldthm1}-(a) a bit further.
Indeed, using that first-layer functions form a vector space,
we obtain a large deviation principle of the variables
$X_N(f)$ jointly in any finite number of $f$'s.
More precisely, for any choice $f_1,\ld, f_k$ first-layer functions, the
random vector $(X_N(f_1),\ld X_N(f_k))$ satifsfies the large deviation principle with rate function
\[
I_{f_1,\ld ,f_k} (x_1,\ld, x_k)= \sup\left(\sum_{i=1}^k x_i t_i- \caP^{\scriptscriptstyle{M}}\left(\sum_{i=1}^k t_i f_i \Big |\mising\right)\right).
\]
We can then take the projective limit, i.e., induce on the space of probability measures $\mathscr{P} (\Omega)$ the topology induced by the maps $\mu\mapsto \int f d\mu$ with $f$ a first-layer function. Then by Dawson-G\"artner theorem \cite[p. 162]{dz}, we have that the random measures
\[
\frac1N\sum_{i=1}^N \delta_{T_i\si}
\]
satisfy the large deviation principle with rate function
\[
\caI (\lambda|\mising)= \sup_{f:\ \text{first-layer function}}\left( \int f d \lambda- \caP^{\scriptscriptstyle{M}}(f|\mising)\right).
\]
This can be considered as the analogue of relative entropy density.

Let us now consider some applications of theorem \ref{ldthm1}.
For the choice $f(\si)= \si_1$ we have $f^* (\tau^1)=\tau^1_0$
\[
\sum_{i=1}^N T_i f= \sum_{i=1}^N \si_i
\]
i.e., we have the large deviation principle and the central limit theorem for the magnetization of the multiplicative Ising model.
Choosing $f(\si)=\si_2$  we have $f^*(\tau^1)=\tau^1_1$, and more generally
$f_k(\si)=\si_{2^k}$,  we have $f^*_k(\tau^1)=\tau^1_k$, i.e., we have the large deviation principle for sums of the form
\[
\frac1N\sum_{i=1}^N \si_{i2^k}
\]
i.e., for the magnetization along decimated lattices.

For the choice $f(\si)= \si_1\si_2$ we have $f^* (\tau^1)=\tau^1_0\tau^1_1$ and
\[
\sum_{i=1}^N T_i f= \sum_{i=1}^N \si_i\si_{2i}.
\]
The function $f(\si)=\si_1\si_3$ is however not a first-layer function and therefore, the large deviations of
\[
\frac1N\sum_{i=1}^N \si_i\si_{3i}
\]
do not follow from Theorem \ref{ldthm1}.

\section{One-dimensional multiplication-invariant Gibbs measures}

The theory developed so far for the multiplicative Ising model quite easily generalizes to one-dimensional
multiplication-invariant Gibbs measures of $\si$ spins, such that the corresponding layers of $\tau$ spins are in the uniqueness regime.
Informally speaking, this means Gibbs measures with formal Hamiltonians
\be\label{formalgenham}
\sum_{i\in \N}\sum_{A}J_A  \prod_{j\in A} \si_{i2^j}
\ee
where the second sum runs over finite subsets of $\N_0$.
E.g. formal Hamiltonians
\[
\sum_{i\in \N} \si_i\si_{2i}+\si_{i}\si_{4i}+\si_{i}\si_{2i}\si_{8i}
\]
are included but not e.g.
\[
\sum_{i\in\N} \si_i\si_{2i}+\si_i\si_{3i}
\]
which will later be called a two-dimensional model.
We choose here to work with powers of $2$ in \eqref{formalgenham}, this can be replaced without any further difficulty by
any prime number. The essential point is that in  \eqref{formalgenham} only powers of a single prime number appear, which makes the
models one-dimensional.

\subsection{One-dimensional potentials}

To define the Gibbs measures with formal Hamiltonian \eqref{formalgenham} more precisely, we define a potential
$U(A,\si)$ to be a function of finite subsets $A$ of $\N$ such that
\ben
\item $U(A,\si)$ depends only on $\si_A$.
\item $\sum_{ A\ni i} \max_{\si_A} |U(A,\si)|$ is finite for all $i\in \N$.
\een
We call such a potential multiplication invariant if
\[
U(iA, \si)=U(A, T_i\si)
\]
for all $A\subset \N$, $\si\in \Omega$, $i\in\N $.
To construct examples of multiplication invariant potentials, we can start from a ``base'' collection
$\{J_A, A\in \caA\}$ of interactions and then define
\be\label{ja}
U(iA, \si)=J_A \prod_{j\in A}T_i\si(j)
\ee
and $U(B,\si)=0$ for $B$ not of the form $iA, A\in \caA, i\in\N$. For the multiplicative Ising model we had
$\caA= \{\{ 1,2\}, \{1\}\}, J_{\{1,2\}}=-\beta J, J_{\{1\}}=-\beta h$. We will from now on restrict to such potentials, which
is the case of $\pm 1$ spins is not a restriction.

We call a potential one-dimensional if the set $\cup_{A\in \caA} A$ contains powers of at most a single prime, which we choose here, without
loss of generality to be $2$, in other words if
\[
\cup_{A\in\caA} A=\{ 2^{i_1}, 2^{i_2}, \ldots\}.
\]
We then have the natural correspondence of the potential $U$ in \eqref{ja} with the shift-invariant potential
\[
V(A+i, \si)= J_A \prod_{j\in A}\theta_i \si(j)
\]
for $i\in\N_0$, $\si\in \Omega_0$,
and $V(B, \si)=0$ for sets not of the form $i+A$.

We call a multiplication-invariant potential $U(A,\si)$ {\em layer unique } if for the corresponding potential $V$ there is a unique
Gibbs measure $\mu^{\scriptscriptstyle{V}}_\infty$ on $\Omega_0$ in the sense of \cite{geo}. Notice that the configuration space $\Omega_0$ corresponds to free boundary conditions at the left end, and so in general despite shift-invariant potentials $V$, the corresponding unique Gibbs measure will  not necessarily be stationary under the shift.

We then have, in complete analogy with theorem \ref{layerthm}, the following result.
\bt\label{genlayerthm}
Let $U$ be a multiplication invariant one-dimensional layer unique potential.
Then $U$ admits a unique Gibbs measure $\mu^{\scriptscriptstyle{U}}_\infty$ which is multiplication invariant.
Under this measure, the layer spins defined by
\[
\tau^r_i = \si_{r2^i}
\]
$r\in 2\N+1$
are independent for different $r$ and distributed
as the unique Gibbs measure with potential $V$ on $\Omega_0$.
\et
As an example, consider $V$ the long-range Ising model:
$V(\{i,j\}, \si)=J(|j-i|) \si_i\si_j$, with $\sum_n nJ(n)<\infty$. Then the corresponding multiplication invariant potential $U$
is given by $U(\{r2^i, r2^j\}, \si)= \si_{r2^i}\si_{r2^j} J(|j-i|)$.

\subsection{Large deviations in the general one-dimensional layer unique context}

The large deviation properties of sums of the form $\frac1N \sum_{i=1}^N T_i f$ are obtained just as in the case of
the multiplicative Ising model. I.e., defining for $f$ a first-layer function the pressure
\[
\caP^M(f|\mu^{\scriptscriptstyle{U}}_\infty)=\lim_{n\to\infty}
\frac1N \log\E_{\mu^{\scriptscriptstyle{U}}_\infty} \left(e^{\sum_{i=1}^N T_i f}\right)
\]
we have
\[
\caP^M(f|\mu^{\scriptscriptstyle{U}}_\infty)=\sum_{k=0}^\infty \frac{1}{2^{k+2}} \; P^{k}_{\mu^{\scriptscriptstyle{V}}_\infty} (f^*)
\]
where
\[
P^{k}_{\mu^{\scriptscriptstyle{V}}_\infty} (f^*)=
\log\E_{\mu^{\scriptscriptstyle{V}}_\infty}\left(e^{\sum_{i=0}^{k} \theta_i f^*}\right).
\]
As a consequence, under $\mu^{\scriptscriptstyle{U}}_\infty$ we have the same results as in theorem
\ref{ldthm1}.

\section{Higher-dimensional models in the uniqueness regime}

Let us start now from a multiplication-invariant potential as constructed in \eqref{ja}
from a collection $\{J_A, A\in\caA\}$. Let us denote, for a such a collection the set
$\cup_{A\in \caA} A= S(J)$. Let us further denote by $P(J)$ the set of primes appearing in the prime factorization of
all the numbers appearing in $S(J)$. We assume that $P(J)$ is a finite set.
We denote by $d=d(J)=|P(J)|$ and call this the dimension of the underlying model
and we order $P(J)= \{p_1,p_2,\ld, p_d\}$ with $p_1<p_2<\cdots<p_d$.
The analogue of the layer decomposition then goes as follows: we write every number
$i\in \N$ in a unique way as
\[
i= r \prod_{i=1}^d p_i^{x_i}
\]
where $x_i\in \N_0$, $r\in \N$, $r$ not divisible by any of the primes $p_1,\ldots, p_d$
(the set of all such $r$ is denoted by $\caK(p_1,\ld,p_d)$).
We further denote, for $N\in \N$ by $\la^r_{p_1,\ld, p_d; N}$ the set
\[
\la^r_{p_1,\ld, p_d; N}= \big\{ (x_1,\ld,x_d)\in (\N_0)^d: r\prod_{i=1}^d p_i^{x_i}\leq N\big\}.
\]
We then have that the Gibbs measure in the lattice interval $[1,N]$ associated to the potential $U$,
with boundary condition $\eta$ factorizes into a product over $r\in \caK(p_1,\ld,p_d), r\leq N$ of independent Gibbs measures
on the sets $\{-1,+1\}^{\la^r_{p_1,\ld, p_d; N}}$ associated to
the corresponding shift-invariant potential $V$, with free boundary conditions on the ``left'' and more ``complicated''
$\eta$-dependent, (and for our purposes here unimportant) boundary conditions on the other ends.

We assume now the following:
\bd
We call the potential $U$ layer unique if the corresponding potential $V$ has a unique infinite-volume Gibbs
measure $\mu^{\scriptscriptstyle{V}}_\infty$ on $\Omega_{0,d}=\{-1,1\}^{(\N_0)^d}$.
\ed
From the layer decomposition of the finite-volume Gibbs measures, and the uniqueness of infinite-volume limits on each layer, we obtain the following analogue of theorem \ref{genlayerthm}.
\bt\label{layeruniquenessgeneralthm}
Let $U$ be a multiplication-invariant potential with associated shift-invariant potential $V$.
Assume that $U$ is layer unique. Then there exists a unique Gibbs measure $\mu^{\scriptscriptstyle{U}}_\infty$
associated to $U$ on the configuration space $\Omega=\{-1,1\}^{\N}$. Under this measure $\mu^{\scriptscriptstyle{U}}_\infty$, the $\tau$ spins defined by
\be\label{kokod}
\tau^{r}_{x_1,\ld,x_d}= \si_{r\prod_{i=1}^d p_i^{x_i}},
\ee
$r\in\caK(p_1,\ld, p_d)$, $x_i\in \N_0$, form independent copies (with respect to $r$) of the measure $\mu^{\scriptscriptstyle{V}}_\infty$.
\et

\subsection{Pressure and large deviations}

In order to obtain the analogue of theorem \ref{ldthm1} in the $d$-dimensional case, we need
the following lemma which is proved in \cite{kiki}.
\bl\label{KIE}
Let $\phi:\N_0\to\R$ be such that there exist $C>0$ and $q>0$ such that $|\phi(n)|< n^q$ for all $x\in\N_0$.
Then there exist a constant $\kappa\in (0,1)$ and functions, $\rho^+,\rho^-:\N\to [0,\infty)$ such that $\rho^+(\ell)>\rho^-(\ell)>0$ for all $\ell$ and such that
\be\label{kieferstuff}
\lim_{N\to\infty}\frac1N \sum_{\substack{r\in \caK(p_1,\ld,p_d)\\ 1\leq r\leq N}} \phi(|\la^r_{p_1,\ld,p_d;N}|)
= \kappa\sum_{j=1}^\infty\left( e^{-\rho^- (j)}- e^{-\rho^+(j)}\right) \phi(j).
\ee
Here $\rho^-,\rho^+$ are defined by
\beq
\nonumber
\rho^-(\ell)&=&\inf \{\rho\geq 0: |D(\rho)|=\ell\}
\nonumber\\
\rho^+(\ell) &=& \sup \{\rho\geq 0: |D(\rho)|=\ell\}
\nonumber
\eeq
with
\[
D(\rho)= \left\{ (x_1,\ld, x_d):\sum_{i=1}^d x_i\log(p_i)\leq \rho\right\}
\]
and
\[
\kappa= \kappa(p_1,\ld,p_d)=1-\frac12-\frac13 +\frac{1}{2\times 3} -\frac15 + \frac{1}{2\times 5} +\frac{1}{3\times5} +\cdots +(-1)^d\frac{1}{p_1\ld p_d}\cdot
\]
\el
\br
Notice that in general we have the bound $\rho^-(\ell)\geq (\ell^{1/d}-1)\log(2)$, ensuring the absolute convergence of the series in
\eqref{kieferstuff}.
In the particular case $d=1$, $p_1=2$, we have $e^{-\rho^- (j)}- e^{-\rho^+(j)}= \frac{1}{2^{j+1}}$, $\kappa=1/2$, consistent with our previous result (lemma \ref{KOROA}).
\er
We now define the analogue of the finite-volume pressures on layers.
For a first-layer function $f:\Omega\to\R$, i.e., a continuous function depending only on $\tau^1$ defined in \eqref{kokod}
($f(\si)= f^*(\tau^1)$), we define
\[
\caP_{\la^r_{p_1,\ld,p_d;N}}(f^*|V)=\E_{\mu^{\scriptscriptstyle{V}}_\infty} \left(e^{\sum_{x\in \la^{r}_{p_1,\ld ,p_d;N}}\theta_x f^*}\right).
\]
Notice that this function, as a function of $r$ and $N$ only depends on $|\la^{r}_{p_1,\ld ,p_d;N}|$ (cf. \cite{kiki}).
Hence, if $|\la^{r}_{p_1,\ld ,p_d;N}|=\ell$, we define
\[
\Psi_\ell(f^*|V)=\caP_{\la^r_{p_1,\ld,p_d;N}}(f^*|V).
\]
We can therefore use lemma \ref{KIE} to obtain the following result.
Define the pressure of a first-layer function w.r.t. the unique Gibbs measure with potential $U$ as before
\be\label{pmd}
\caP^{\scriptscriptstyle{M}}(f|\mu^{\scriptscriptstyle{U}}_\infty)=
\lim_{N\to\infty} \frac1N \log \E_{\mu^{\scriptscriptstyle{U}}_\infty} \left(e^{\sum_{i=1}^N T_i f}\right).
\ee
\bt\label{oula}
Let $U$ be a layer-unique $d$-dimensional multiplication invariant potential, and $f$ a first-layer function. Then the limit defining the pressure \eqref{pmd} exists and is given by
\[
\caP^{\scriptscriptstyle{M}}(f|\mu^{\scriptscriptstyle{U}}_\infty)=
\kappa\sum_{j=1}^\infty\left( e^{-\rho^- (j)}- e^{-\rho^+(j)}\right) \Psi_{j} (f^*|V).
\]
As a consequence, the random variables
\[
X_N (f)= \frac1N \sum_{i=1}^N T_i f
\]
satisfy the large deviation principle with rate function
\[
I_f(x)=\sup_{t\in\R}\left(tx-\caP^M(t f|\mu^{\scriptscriptstyle{U}}_\infty)\right).
\]
\et
\br
Just as in the one-dimensional case, by considering joint large deviations of $X_N(f_1),\ld,X_N(f_k)$ and taking a projective limit, we have the large deviation principle for the random probability measures
\[
\frac1N\sum_{i=1}^N \delta_{T_i\si}
\]
in the weak topology induced by first-layer functions.
\er

\subsection{Dimensional extension and general large deviations}

We can extend the large deviation principle in the following way. Suppose e.g. we want to obtain the large deviation principle
of
\be\label{qwe}
\frac1N\sum_{i=1}^N \si_i\si_{2i}+\si_i\si_{3i}
\ee
under the multiplicative Ising model with $h=0, \beta J=1$ (for simplicity). We can view this multiplicative Ising model
as the model with formal Hamiltonian
\[
\sum_i \si_{i} \si_{2i} + 0.\si_i \si_{3i}
\]
i.e., a two-dimensional model consisting in each layer of independent copies of $\mising$.
The corresponding two-dimensional shift-invariant potential is
\[
V(\{(x_1,x_2), (x_1+1, x_2)\}, \tau)= \tau_{x_1}\tau_{x_1+1}
\]
and $V(A,\tau)=0$ for other subsets $A$,
i.e., the potential has only interaction in the $x_1$ direction.
This model is of course still
layer unique, and hence we have the large deviation principle for \eqref{qwe}, because it is now a first-layer function in the two-dimensional model.

More generally, suppose that we have a layer unique Gibbs measure $\mu^{\scriptscriptstyle{U}}_\infty$ corresponding to a $d$ dimensional multiplication invariant potential,
and we want to prove the large deviation principle for
\[
\frac1N\sum_{i=1}^N  T_i f
\]
where $f$ is a local function, i.e. a function depending only on a finite number of coordinates.
Since $f$ is local we can write
\[
f= \sum_{B} J_B \si_B
\]
where $\si_B= \prod_{i\in B}\si_i$ and where the sum over $B$ runs over a finite number of finite subsets $B\subset \N$.
Let us call $\pi_1, \ld, \pi_k$ the primes involved in the prime decomposition of $\cup B$ and $p_1,\ld, p_d$ the primes
associated to the potential $U$. Then, to study the large deviation of $\frac1N\sum_{i=1}^N T_i f$, we go to a higher dimensional model
associated to the primes $\{\pi_1, \ld, \pi_k\}\cup \{p_1,\ld, p_d\}:= \{p_1,\ld, p_d, p'_{d+1}, \ld, p'_{d+d'}\}$, consisting of
$d'$ non-interacting layers distributed according to $\mu^{\scriptscriptstyle{U}}_\infty$. This new model is of course still
layer-unique (as we have not added any new interaction), and $f$ is a first-layer function in the new model. Therefore,
$\frac1N\sum_{i=1}^N  T_i f$ satisfies the large deviation principle under $\mu^{\scriptscriptstyle{U}}_\infty$.

As a conclusion, for a layer-unique multiplication-invariant potential $U$, we have the large deviation principle for
$(\frac1N\sum_{i=1}^N  T_i f)_N$ for every local $f$. 
Because the set of local functions is dense in the set of continuous function for the sup-norm topology, we can 
summarize our observations in the following theorem.
\bt
Let $\mu^{\scriptscriptstyle{U}}_\infty$ be a layer-unique Gibbs measure with a multiplication-invariant potential $U$.
For any local function $f$, $\caP^{\scriptscriptstyle{M}}(f|\mu^{\scriptscriptstyle{U}}_\infty)$ exists and the sequence
$(\frac1N\sum_{i=1}^N T_i f)_N$ satisfies the large deviation principle with rate function
\[I_f(x)=\sup_{t\in\R}\left(tx-\caP^M(t f|\mu^{\scriptscriptstyle{U}}_\infty)\right).
\]
As a consequence, by taking the projective limit, we have that the sequence of random measures
$(\frac1N\sum_{i=1}^N \delta_{T_i\si})_N$ satisfies the large deviation principle in the weak topology with rate function
\[
\caI (\lambda|\mising)= \sup_{f\ \textup{local}}\left( \int f d \lambda- \caP^{\scriptscriptstyle{M}}(f|\mu^{\scriptscriptstyle{U}}_\infty)\right).
\]
\et

\subsection{A Shannon-McMillan-Breiman theorem}

In Subsection \ref{KS} we computed the Kolmogorov-Sinai entropy of the multiplicative Ising model.
Here we prove the corresponding almost-sure convergence in the more general context of
layer-unique Gibbs measures in dimension one. Then we explain how to
extend it, in the same spirit as we did for our large deviation results.

\bt
Let $U$ be a layer-unique one-dimensional multiplication invariant potential and
$\mu^{\scriptscriptstyle{U}}_\infty$ the corresponding Gibbs measure.
Then $\mu^{\scriptscriptstyle{U}}_\infty$-almost surely, we have
\be\label{ouou}
\lim_{N\to\infty}- \frac1N \log \mu^{\scriptscriptstyle{U}}_\infty(\si_{[1,N]})=
\sum_{k=0}^\infty \frac{1}{2^{k+2}}\ s^{\scriptscriptstyle{V}}_{k+1}
\ee
where
\[
s^{\scriptscriptstyle{V}}_{k+1}=-\E_{\mu_\infty^{\scriptscriptstyle{V}}}
\left(\log\mu_\infty^{\scriptscriptstyle{V}}(\tau_0,\ldots,\tau_k)\right).
\]
\et
\bpr
We are going to use the large deviation principle for the sequence
$(- \frac1N \log \mu^{\scriptscriptstyle{U}}_\infty(\si_{[1,N]}))_N$ and
consequent strong law of large numbers.\\
We again use the following formula
\[
\log \mu^{\scriptscriptstyle{U}}_\infty (\si_{[1,N]})=
\sum_{\substack{1\leq r\leq N\\  r \ \text{odd}}}\log\mu^{\scriptscriptstyle{V}}_\infty
(\tau^r_0,\ldots, \tau^r_{\psi_2 (r/N)}).
\]
As a consequence, using the fact that different layers are independent and have the same distribution
$\mu^{\scriptscriptstyle{V}}_\infty$, we have for every $t\in\R$
\begin{align*}
F_N(t)
& := \log \E_{\mu_\infty^{\scriptscriptstyle{U}}}
\left( e^{-t \log \mu^{\scriptscriptstyle{U}}_\infty(\si_{[1,N]})}\right)\\
& = \sum_{\substack{1\leq r\leq N\\  r \ \text{odd}}} \log 
\E_{\mu_\infty^{\scriptscriptstyle{V}}}
\left( e^{-t \log \mu^{\scriptscriptstyle{V}}_\infty(\tau_0,\ld, \tau_{\psi_2 (r/N)})}\right).
\end{align*}
Because $\mu^{\scriptscriptstyle{V}}_\infty$ is a Gibbs measure, there exist strictly positive
constants $c_1,c_2$ such that for all $k$
\be\label{Gigi}
e^{-c_1 k} \leq \mu^{\scriptscriptstyle{V}}_\infty(\tau_0,\ld, \tau_{k})\leq e^{-c_2 k}.
\ee
It follows from lemma \ref{KOROA} that 
\[
F(t):=\lim_{N\to\infty} \frac{1}{N} F_N(t)=
\sum_{k=0}^\infty \frac{1}{2^{k+2}}\ s^{\scriptscriptstyle{V}}_{k+1}(t)
\]
where
\[
s^{\scriptscriptstyle{V}}_{k+1}(t):=
\log \E_{\mu^{\scriptscriptstyle{V}}_\infty}
\left[ (\mu^{\scriptscriptstyle{V}}_\infty(\tau_0,\ld, \tau_{k}))^{-t}\right].
\]
The map $t\mapsto s^{\scriptscriptstyle{V}}_{k+1}(t)$ is continuously differentiable and
strictly convex for every $k$.
Moreover,
\begin{align*}
\left|\frac{\textup{d}}{\textup{d} t} s^{\scriptscriptstyle{V}}_{k+1}(t) \right|
&= \left|\frac{\E_{\mu^{\scriptscriptstyle{V}}_\infty}\left(e^{t Y_k} Y_k \right)}{
\E_{\mu^{\scriptscriptstyle{V}}_\infty}\left(e^{t Y_k} \right)}\right| \\
& \leq c_1 k
\end{align*}
where $Y_k=-\log \mu^{\scriptscriptstyle{V}}_\infty(\tau_0,\ld, \tau_{k})$. Therefore $F$ is
strictly convex and continuously differentiable. 
By G\"artner-Ellis theorem (see \cite{dz}) the sequence
$(- \frac1N \log \mu^{\scriptscriptstyle{U}}_\infty(\si_{[1,N]}))_N$ satisfies a large deviation
principle with a strictly convex rate function. As a consequence, we have exponential convergence
and by the Borel-Cantelli lemma, the strong law of large numbers.
\epr

\br\leavevmode\\
\begin{enumerate}
\item 
In the spirit of Theorem \ref{oula}, we can easily extend the previous theorem to layer-unique $d$-dimensional multiplication-invariant potentials, but the analogue of formula \eqref{ouou} is quite
cumbersome.
\item
The theorem is valid beyond Gibbs measures, namely if the different layers are independent and have the same distribution and estimates \eqref{Gigi} hold, then the same proof works.
\end{enumerate}
\er

\bigskip

\noindent\textbf{Acknowledgments.} 
This work was completed during a 3-month stay of F. R. at the Centre de Physique Th\'eorique
funded by the CNRS. J.-R. C. thanks Mike Hochman for enlightening dicussions and Benjy Weiss for 
having pointed to him reference \cite{jen}.


\end{document}